\documentclass[12pt, reqno]{amsart}
\usepackage[dvipsnames]{xcolor}
\usepackage{mathtools}
\mathtoolsset{showonlyrefs}
\usepackage{etoolbox}
\usepackage{blindtext}
\usepackage{latexsym}
\usepackage{txfonts}
\usepackage{lmodern}
\usepackage[pagewise]{lineno}
\usepackage[utf8]{inputenc}
\usepackage{exscale}
\usepackage{amsmath}
\usepackage{amssymb}
\usepackage{amsfonts}
\usepackage{mathrsfs}
\usepackage{amsbsy}
\usepackage{relsize}
\usepackage{comment}
\PassOptionsToPackage{reqno}{amsmath}
\usepackage{esint} % For \fint symbol
\usepackage{amssymb} % For \lesssim symbol
\usepackage{stmaryrd}
\usepackage{cite} 

\usepackage{tikz}
\usepackage[top=2.5cm, bottom=2.5cm, left=2.5cm, right=2.5cm]{geometry}
\usepackage{mathtools}
\usepackage{soul}
\usepackage{color}

\usepackage{slashed}
\definecolor{mypink1}{rgb}{0.858, 0.188, 0.478}
\definecolor{mypink2}{RGB}{219, 48, 122}
\definecolor{mypink3}{cmyk}{0, 0.7808, 0.4429, 0.1412}
\definecolor{mygray}{gray}{0.6}
\definecolor{venetianred}{rgb}{0.78, 0.03, 0.08}
\definecolor{sapphire}{rgb}{0.03, 0.15, 0.4}
\definecolor{utahcrimson}{rgb}{0.83, 0.0, 0.25}
\definecolor{trueblue}{rgb}{0.0, 0.45, 0.81}
\definecolor{carminered}{rgb}{1.0, 0.0, 0.22}
\definecolor{cobalt}{rgb}{0.0, 0.28, 0.67}
\definecolor{cornflowerblue}{rgb}{0.39, 0.58, 0.93}
\definecolor{darkmagenta}{rgb}{0.55, 0.0, 0.55}
\definecolor{electricultramarine}{rgb}{0.25, 0.0, 1.0}
\definecolor{falured}{rgb}{0.5, 0.09, 0.09}
\definecolor{hancornflowerblue}{rgb}{0.32, 0.09, 0.98}
\definecolor{mahogany}{rgb}{0.75, 0.25, 0.0}
\definecolor{oucrimsonred}{rgb}{0.6, 0.0, 0.0}
\definecolor{persianblue}{rgb}{0.11, 0.22, 0.73}
\definecolor{rufous}{rgb}{0.66, 0.11, 0.03}
\definecolor{uablue}{rgb}{0.0, 0.2, 0.67}
\definecolor{zaffre}{rgb}{0.0, 0.08, 0.66}
\definecolor{carmine}{rgb}{0.59, 0.0, 0.09}
\usepackage[colorlinks, linkcolor=carminered, citecolor=electricultramarine, urlcolor=mahogany, hypertexnames=false]{hyperref}
\newcommand{\R}{\ensuremath{\mathbb{R}}}

\newcommand{\N}{\ensuremath{\mathbb{N}}}

\newtheorem{definition}{Definition}[section]
\newtheorem{proposition}{Proposition}[section]
\newtheorem{theorem}{Theorem}[section]
\newtheorem{remark}{Remark}[section]
\newtheorem{lemma}{Lemma}[section]

\title[Stochastic parabolic equation]{The inhomogeneous fractional stochastic heat equation driven by fractional Brownian motion}

\author[R. Alessa, R. Al Subaie, M. Alwohaibi, M. Majdoub, E. Mliki]{\tiny Rasha Alessa, Reem Al Subaie, Maram Alwohaibi\\\tiny Mohamed Majdoub, Ezzedine Mliki}
\address{\small Department of Mathematics, College of Science, Imam Abdulrahman Bin Faisal University\\ P. O. Box 1982, Dammam, Saudi Arabia}
\address{\small Basic and Applied Scientific Research Center, Imam Abdulrahman Bin Faisal University\\ P.O. Box 1982, 31441, Dammam, Saudi Arabia}
\email{\sl \color{blue}{ralessa@iau.edu.sa}}
\email{\sl \color{blue}{rmalsubaei@iau.edu.sa}}
\email{\sl \color{blue}{malwohaibi@iau.edu.sa}}
\email{\sl {\color{blue}{mmajdoub@iau.edu.sa}} ({\tt Corresponding author})}
\email{\sl \color{blue}{ermliki@iau.edu.sa}}

\begin{document}
\begin{abstract} We investigate  the fractional Hardy-H\'enon equation with fractional Brownian noise
$$
\partial_tu(t)+(-\Delta)^{\theta/2} u(t)=|x|^{-\gamma} |u(t)|^{p-1}u(t)+\mu \, \partial_t B^H(t),
$$
where $\theta>0$, $p>1$, $\gamma\geq 0$, $\mu \in\R$, and the random forcing $B^H$ is the fractional Brownian motion defined on some complete probability space $(\Omega, \mathcal{F}, \mathbb{P})$ with Hurst parameter $H\in (0,1)$.
We establish the local existence and uniqueness of mild solutions under appropriate conditions on the parameters of the equation. 
\end{abstract}

%@@@@@@@@@@@@@@@@@@@@@@@@@@@@@@@@@@@@%@@@@@@@@@@@@@@@@@@@@@@@@@@@@@@@@@@@@


%@@@@@@@@@@@@@@@@@@@@@@@@@@@@@@@@@@@@%@@@@@@@@@@@@@@@@@@@@@@@@@@@@@@@@@@@@

%\date{\today}

\maketitle
%%%%%%%%%%%%%%%%%%%%%%%%%%%%%%%%%%%%%%%%%%%%%%%%%%%%%%%%%%%%%%%%%%%%%%%%%%%%%%%%%%%%%%%%%%%%%%%Ã¹%%%

\section{Introduction}
\label{S1}

We deal with the initial value problem for the fractional Hardy-H\'enon equation with fractional Brownian noise
\begin{equation}
\left\{ \begin{array}{ll}
     \label{main}
\partial_tu(t)+(-\Delta)^{\theta/2} u(t)=|x|^{-\gamma} |u(t)|^{p-1}u(t)+\mu\, \partial_t B^H(t), \\
    u(0) = u_0, 
\end{array}
\right.
\end{equation}
where $\theta>0$, $p>1$, $\gamma\geq 0$, $\mu \in\R$, and the random forcing $B^H$ is the fractional Brownian motion defined on some complete probability space $(\Omega, \mathcal{F}, \mathbb{P})$ with Hurst parameter $H\in (0,1)$. When $H = 1/2$, $B^{1/2}$ corresponds to the standard Brownian motion.

As is customary, we analyze \eqref{main} by considering the associated integral equation:
\begin{equation}
\label{integral}
\begin{split}
u(t)= {\rm e}^{-t(-\Delta)^{\theta/2}}u_{0}&+\int_{0}^{t}{\rm e}^{{-}(t-s)(-\Delta)^{\theta/2}}\,\left(|x|^{-\gamma}|u(s)|^{p-1}u(s)\right)\,ds\\&+\mu\int_0^t{\rm e}^{{-}(t-s)(-\Delta)^{\theta/2}}\,dB^H(s),
\end{split}
\end{equation}
where ${\rm e}^{-t(-\Delta)^{\theta/2}}$ is the linear heat semi-group.
%%%%%%%%%%%%%%%%%%%%%%%%%%%%%%%%%%%%%%%%%%%%%%%%%%%%%%%%%%%%%%%%%%%%%%%%%%%%%%%%%%%%%%%%%%%%%%%Ã¹%%%%%%%%

The initial value problem \eqref{main} without fractional noise (i.e., when $\mu=0$):
\begin{equation}
\label{mu=0}
\partial_t{u}+(-\Delta)^{\theta/2} {u}=|x|^{-\gamma} |u|^{p-1}u,
\end{equation}
has garnered significant attention within the mathematical community. In particular, the case $\theta=2$ and $\gamma=0$ corresponds to the classical semi-linear heat equation, which has been extensively studied (see, for example, \cite{Br, HW, MI, QS, Rib, W1, W2, W3}). For nonlinearities that grow faster than a power, further investigations can be found in \cite{IJMS, Ioku, IRT, MT, MT1}. When considering $\theta=2$ and $\gamma>0$, known as the Hardy-H\'enon parabolic model, we recommend referring to \cite{BTW, Tay} and the references provided therein for more detailed discussions and analysis.

In \cite{MY}, the global well-posedness of solutions to \eqref{mu=0} for general values of $\theta$ was extensively examined. The study primarily focused on small initial data within pseudomeasure spaces. Furthermore, an unified method was introduced in \cite{MYZ2008} to handle the well-posedness of \eqref{mu=0} in the case of $\gamma=0$ for general nonlinearities. Note that the time-fractional scenario with time-space white noise has been studied in \cite{Hach}.

It is worth noting that the case $\gamma>0$ for $\theta=2$ was studied in \cite{BTW}. In our analysis, we will make use of the approach employed in \cite{BTW}, along with the properties of the fractional heat kernel as described in \cite{MYZ2008}. By combining these techniques, we aim to address our specific problem and gain valuable insights into its solution.

The development of stochastic calculus for fractional Brownian motion ({\sf fBm}) has naturally paved the way for investigating stochastic partial differential equations (SPDEs) driven by {\sf fBm}. The study of such SPDEs represents a significant research avenue within probability theory and stochastic analysis, and it has yielded numerous insightful findings. This line of research has stemmed from the extensive applications of {\sf fBm} in various fields, motivating the exploration of SPDEs as a powerful tool to model and analyze complex systems.

In recent years, the investigation of SPDEs driven by fractional noise has garnered significant attention within the mathematical community. This growing interest stems from both theoretical considerations and the broad range of applications in fields such as physics, biology, hydrology, and other scientific disciplines. Of particular significance is the study of the well-posedness of semilinear stochastic parabolic equations driven by infinite-dimensional fractional noise, which has been a subject of special interest. Noteworthy contributions in this area can be found in works such as \cite{DDM, MS, ND, SV, Zhang}.

Moreover, the exploration of alternative types of noise has also been undertaken in the literature, as evidenced by studies conducted in \cite{BT, BNS, CM, HN09}. 
{Furthermore, time-space fractional stochastic parabolic equations driven by fractional noise have been studied in \cite{Duan}. More recently, a novel stochastic fractional partial differential equation was introduced in \cite{Zili}. This equation features a mixed operator combining the standard Laplacian, the fractional Laplacian, and the gradient operator. It is driven by a random noise with a covariance measure structure in time and behaves as a Wiener process in space.}

These investigations further enhance our understanding of the dynamics and behavior of SPDEs, offering valuable insights into their mathematical properties and practical applications.

In this context, we focus on studying the semilinear Hardy-Hénon equation driven by fractional Brownian noise given by \eqref{main}. Our primary objective is to find local-in-time solutions to \eqref{main} for initial data belonging to a specific Lebesgue space. To achieve this, we introduce the concept of a mild solution for \eqref{main}, which provides a suitable framework for characterizing the solution behavior. {Note that the solutions \( u(t) \) are random fields, meaning they depend on both space and time, as well as the random parameter \( \omega \in \Omega \). Consequently, \( u(t) \) is a random variable taking values in a function space. The Lebesgue norms \(\| \cdot \|_q\), for \( 1 \leq q \leq \infty \), are spatial norms applied to the spatial component of \( u(t) \) for each fixed \( t \) and \( \omega \).}
\begin{definition}
\label{mild}
Let $T>0$ and $1\leq q\leq \infty$. A measurable function $u : \Omega\times [0,T]\to L^q$ is a mild solution of \eqref{main} if
\begin{itemize}
\item[i)] $u\in C([0,T]; L^q)$,
\item[ii)] $u$ satisfies \eqref{integral} with probability one.
\end{itemize}
\end{definition}

We have successfully established the following result regarding the local well-posedness of the problem at hand.
\begin{theorem}
\label{LWP}
Let $N\geq 1$, $p>1$, and $\mu\neq 0$. Assume that  
\begin{equation}
    \label{Main-Ass}
    \begin{split}
        2\theta>1,\quad   &0 < \gamma < \min(\theta,\, N),\\
        \max\bigg(\frac{1}{q},\, \frac{N}{2\theta},\, \frac{1}{2} \bigg)&<H<1, \quad\max\bigg(\frac{Np}{N-\gamma },\, \frac{N(p-1)}{\theta-\gamma}\bigg)<q<\infty.
    \end{split}
\end{equation}

Then, for any $u_0\in L^q$,  there exists ${T=T(\|{u_0}\|_q, H)>0}$ such that problem \eqref{main} possesses a unique mild solution on $[0, T].$
\end{theorem}
\begin{remark}
~{\rm
\begin{itemize}
    \item[(i)] The Cauchy problem \eqref{main} was studied in \cite{CO} for $\theta=2$, $\gamma=0$ and $\mu=1$. The local existence was proved under restrictive assumptions.
    \item[(ii)] In \cite{MMEM}, the authors improved the results of \cite{CO} by relaxing the assumptions on $\gamma, q$.
    \item[(iii)] The assumption $H> \frac{N}{2\theta}$ ensures that the solution is a space-time continuous random field. More precisely, the parabolic equation is formulated as a stochastic Cauchy problem within a suitable Lebesgue space, and its solution is sought in the mild form, specifically as a stochastic convolution integral. See \cite{JFA2022} for further explanations.
\item[(iv)] { While the Hurst parameter \(H\) does not explicitly appear in the contraction estimates \eqref{Contra-1}--\eqref{Contraction} below, it plays an implicit role in determining the choice of parameters. For instance, the condition \(H > \frac{N}{2\theta}\) in Theorem \ref{LWP} is essential to ensure the space-time continuity of the solution \(Z(t)\), a key requirement for the well-posedness of the problem. Furthermore, the stochastic convolution term \(Z(t)\) inherently depends on \(H\), and its regularity properties significantly influence the behavior of the solution \(u(t)\). Thus, \(H\) indirectly shapes the analysis and the resulting solution framework.}
\item[(v)]{Building on the approach developed in \cite{Tay24, Chik22}, we believe that the use of  weighted spaces could be a suitable strategy to address the singularity in the space variable of the drift term.}
\end{itemize}
}
\end{remark}
\begin{remark}
~{\rm {The main difficulty associated with the {\sf fBm}  in the context of the stochastic heat equation \eqref{main} is indeed the interpretation of the stochastic integral with respect to the {\sf fBm}. This issue arises because {\sf fBm} is not a semimartingale for $H\neq 1/2$, which complicates the definition of the stochastic integral using standard It\^o calculus.} { We address this challenge by employing the properties of {\sf fBm} and the associated stochastic calculus. Specifically, we provide a detailed discussion of the {\sf fBm} and its covariance structure, which is crucial for defining the stochastic integral. We also introduce the Wiener integral representation of {\sf fBm} (see \eqref{BH} below) and discuss the kernel that allows for the construction of the stochastic integral.}}
\end{remark}
The paper is organized as follows. In Section \ref{S2}, we offer a comprehensive review of the essential background information necessary for understanding the subsequent proofs. Specifically, we provide a detailed overview of the relevant concepts and mathematical tools that form the foundation of our analysis. Section \ref{S3} is devoted to the proof of our main result, namely Theorem \ref{LWP}.

{Throughout this article, the symbol $C$ denotes generic positive constants that are not essential to the analysis and may change from line to line.}

\section{Preliminaries}
\label{S2}
\subsection{Fractional Brownian motion}
In the subsequent discussion, $(\Omega, \mathcal{F}, \mathbb{P})$ denotes a complete probability space. Fractional Brownian motion ({\sf fBm}) was initially introduced and investigated by Kolmogorov \cite{Kol} within the framework of Hilbert spaces. Let $[0, T]$ represent a time interval with an arbitrary fixed $T>0$. A {\sf fBm} with a Hurst parameter $H\in (0, 1)$ is a centered Gaussian process denoted as $B^{H}$, possessing the following covariance structure:
\begin{equation}\label{cov}
R(s, t):=E(B^{H}(t)B^{H}(s))= \frac{1}{2} \left(t^{2H}+s^{2H}- |t-s|^{2H}\right),
\end{equation}
where $s,t\in [0, T].$ It is worth emphasizing that when the Hurst parameter takes the value of $H = 1/2$, the process $B^{1/2}(t)$ corresponds to the standard Brownian motion.

The  {\sf fBm} can be alternatively defined as the unique self-similar Gaussian process exhibiting stationary increments
\begin{equation*}
E[(B^{H}(t)-B^{H}(s))^{2}]=|t-s|^{2H},
\end{equation*}
and $H$-self similar
\begin{equation*}
\left(\frac{1}{c^{H}}B^{H}(ct), \, t\geq0\right)\overset{d}{=}\left(B^{H}(t), \, t\geq0\right),
\end{equation*}
for all $c>0.$  Here, the symbol $\overset{d}{=}$ denotes equality in the distributional sense. Furthermore, the process $B^{H}$ can be expressed using the following Wiener integral representation:
\begin{equation}\label{BH}
B^{H}(t)=\int_{0}^{t}K^{H}(t, s)\, dW(s).
\end{equation}
In the above representation, $W={W(t);, t\in [0, T]}$ denotes a Wiener process, and $K^{H}(t, s)$ represents the kernel defined as follows:
 \begin{equation}\label{KH}
 K^{H}(t, s)=c_{H}(t-s)^{H-\frac{1}{2}}+c_{H} (\frac{1}{2}-H)\int_{s}^{t}(r-s)^{H-\frac{3}{2}}(1-(1+(\frac{s}{r}))^{\frac{1}{2}-H})\, dr,
 \end{equation}
 where $c_{H}$ is given by
 \begin{equation*}
 c_{H}=\left(\frac{2H\Gamma(\frac{3}{2}-H)}{\Gamma(H+\frac{1}{2})\Gamma(2-2H)}\right)^{\frac{1}{2}}.
   \end{equation*}
 From \eqref{KH} we obtain
  \begin{equation*}
  \frac{\partial K^{H}}{\partial t}(t, s)=c_{H}\left(H-\frac{1}{2}\right)\left(\frac{s}{t}\right)^{\frac{1}{2}-H}(t-s)^{H-\frac{3}{2}}.
   \end{equation*}
   It is noteworthy that when $H>1/2$, the kernel $K^{H}(t,s)$ exhibits regularity and can be expressed in a simpler form as follows:
   \begin{equation*}
   K^{H}(t,s)=c_{H} s^{\frac{1}{2}-H}\int_{s}^{t}(r-s)^{H-\frac{3}{2}}r^{H-\frac{1}{2}}\,dr.
   \end{equation*}
Let $\varepsilon_{H}$ denote the linear space of step functions on $[0, T]$ with the following form:
 \begin{equation}\label{pphi}
 \varphi(t)=\sum_{i=1}^{n}a_{i}1_{(t_{i}, t_{i+1}]}(t),
  \end{equation}
  where $t_{1}, ... , t_{n}\in [0,T],$ $n\in \mathbb{N},$ $a_{i}\in \mathbb{R}$ and by $\mathcal{H}$ the closure of $\varepsilon_{H}$ with respect to the scalar product
   \begin{equation*}
\langle1_{[0, t]}, \, 1_{[0, s]}\rangle_{\mathcal{H}}=R(t,s).
  \end{equation*}
  For $\varphi \in\varepsilon_{H}$, characterized by the form \eqref{pphi}, we define its Wiener integral with respect to the {\sf fBm} as follows:
  \begin{equation*}
  \int_{0}^{T}\varphi_{s}\,dB^{H}(s)=\sum_{ i=1}^{n}a_{i}(B^{H}(t_{i+1})-B^{H}(t_{i})).
  \end{equation*}
  Evidently, the mapping
\begin{equation*}
\varphi=\sum_{i=1}^{n}a_{i}1_{(t_{i}, t_{i+1}]}\longrightarrow \int_{0}^{T}\varphi_{s}\,dB^{H}(s),
\end{equation*}
serves as an isometry between $\varepsilon_{H}$ and the linear space $\text{span}\{B^{H}(t),\; t\geq 0\}$ when regarded as a subspace of $L^{2}(\Omega)$. The image of an element $\Phi \in \mathcal{H}$ under this isometry is referred to as the Wiener integral of $\Phi$ with respect to $B^{H}$.

For any $s<T$, let us consider the operator $L^{\ast}_H: \mathcal{H}\to L^{2}([0, T])$ given by
          \begin{equation*}
          (L^{\ast}_{H}\varphi)(s)=K^{H}(T,s)\varphi(s)+\int_{s}^{T}(\varphi(r)-\varphi(s))\frac{\partial K^{H}}{\partial r}(r,s)\, dr.
          \end{equation*}
      When $H>1/2$, the operator $L^{\ast}_{H}$ has the simpler expression
      \begin{equation*}
                (L^{\ast}_{H}\varphi)(s)=\int_{s}^{T}\varphi(r)\frac{\partial K^{H}}{\partial r}(r,s)\, dr.
          \end{equation*}
    We refer to \cite{AMN} for the proof that $L^{\ast}_{H}$ constitutes an isometry between $\mathcal{H}$ and $L^{2}([0, T])$. Consequently, we can establish the following relationship between the Wiener process $W$ and the fractional Brownian motion $B^{H}$:
     \begin{equation*}
     \int_{0}^{t}\varphi(s)\,dB^{H}(s)=\int_{0}^{t}(L^{\ast}_{H}\varphi)(s)\, dW(s),
     \end{equation*}
     for every $t\in [0, T]$ and $\varphi 1_{[0, t]}\in \mathcal{H}$ if and only if $L^{\ast}_{H}\varphi \in L^{2}([0, T]).$ Furthermore, it is worth recalling that for $\phi, \chi \in \mathcal{H}$ satisfying $\displaystyle\int_{0}^{T}\int_{0}^{T}|\phi(s)||\chi(t)||t-s|^{2H-2}\, ds dt<\infty$, their scalar product in $\mathcal{H}$ can be expressed as:
\begin{equation}\label{intW}
\langle\phi, \chi\rangle_{H}= H(2H-1)\int_{0}^{T}\int_{0}^{T}\phi(s)\chi(t)|t-s|^{2H-2}\, ds dt.
\end{equation}
This formula establishes the inner product relationship between $\phi$ and $\chi$ in the space $\mathcal{H}$.

       It is important to note that in general, the existence of the right-hand side of equation \eqref{intW} requires careful justification (see \cite{ND}). However, since we will be exclusively working with Wiener integrals over Hilbert spaces, it is worth highlighting that if $X$ is a Hilbert space and $u\in L^{2}([0, T]; X)$ is a deterministic function, the relation \eqref{intW} holds. Moreover, the right-hand side of \eqref{intW} is well-defined in $L^{2}(\Omega, X)$ if $L^{\ast}_{H}u$ belongs to $L^{2}([0, T]\times X)$.
       
          \subsection{Cylindrical fractional Brownian motion}
          Following the approach in \cite{DDM}, we define the standard cylindrical fractional Brownian motion in $X$ as the formal series:
\begin{equation}\label{seri}
B^{H}(t)=\sum_{n=1}^{\infty}e_{n}b_{n}^{H}(t),
\end{equation}
where $\{e_{n},\, n\in \mathbb{N}\}$ represents a complete orthonormal basis in $X$, and $b_{n}^{H}$ denotes a one-dimensional {\sf fBm}.

 It is widely recognized that the infinite series \eqref{seri} does not converge in $L^{2}(\mathbb{P})$, which implies that $B^{H}(t)$ is not a well-defined $X$-valued random variable. However, for any Hilbert space ${\mathcal N}$ such that $X\hookrightarrow {\mathcal N}$, where the embedding of $X$ into ${\mathcal N}$ is a Hilbert-Schmidt operator, the series \eqref{seri} defines a ${\mathcal N}$-valued random variable. Consequently, $\{B^{H}(t), \, t\geq0\}$ can be regarded as a ${\mathcal N}$-valued identity fractional Brownian motion (Id-{\sf fBm}).
 
          Following the approach introduced for cylindrical Brownian motion in \cite{DZ}, it is possible to define a stochastic integral of the form:
\begin{equation}\label{stint}
\int_{0} ^{T}f(t)\,{d}B^{H}(t),
\end{equation}
where $f:[0, T]\mapsto \mathcal{L}(X, Y)$ and $Y$ represents another real and separable Hilbert space. The expression \eqref{stint} corresponds to a random variable valued with $Y$ that remains independent of the choice of ${\mathcal N}$.

          Consider a deterministic function $f$ with values in $\mathcal{L}_{2}(X, Y)$, which represents the space of Hilbert-Schmidt operators from $X$ to $Y$. We now introduce the following assumptions on $f$:
          \begin{itemize}
\item[i)] For each $x\in X$, $f(.)x\in L^{p}([0, T]; Y),$ for $q H>1.$
\item[ii)] $\displaystyle\int_{0}^{T}\int_{0}^{T}\|f(s)\|_{\mathcal{L}_{2}(X, Y)}\,\|f(t)\|_{\mathcal{L}_{2}(X, Y)}\,|s-t|^{2H-2}\,ds dt<\infty.$
\end{itemize}
{From \cite[(2.15), p. 230]{DDM}, it is established that the stochastic integral \eqref{stint} can be properly defined as:}
\begin{equation}\label{sumint}
\int_{0} ^{T}f(t)\,{d}B^{H}(t):=\sum_{n=1}^{\infty}\int_{0}^{{T}} f(t) e_{n}\, db_{n}^{H}(t)=\sum_{n=1}^{\infty}\int_{0}^{{T}}(L_{H}^{\ast}fe_{n})\, db_{n}^{{H}}(t),
\end{equation}
where $b_{n}^{H}$ represents the standard Brownian motion associated with the fractional Brownian motion $b_{n}^{H}$ through the representation formula \eqref{BH}. It is important to note that since $fe_{n}\in L^{2}([0, T]; Y)$ for each $n\in \mathbb{N}$, the variables ${\int_{0}^{T}fe_{n}\, db_{n}^{H}}$ are mutually independent (see \cite{DDM}).
The series \eqref{sumint} is finite and can be expressed as:
          \begin{equation}\label{fiseri}
      \sum_{n=1}^\infty\|L^{\ast}_{H} (fe_{n})\|^{2}=\sum_{n=1}^\infty\parallel\parallel fe_{n}\parallel_{\mathcal{H}}\parallel_{Y}^{2}<\infty.
          \end{equation}
         In the case where $X=Y=\mathcal{H}$, we have:
             \begin{eqnarray*}
         \sum_{n=1}^{\infty}\int_{0} ^{t}f(s)e_{n}db^{H}_{n}(s)&=&\sum_{n=1}^{\infty}\sum_{m=1}^{\infty}e_{m}\int_{0}^{t}\langle f(s)e_{n}, e_{m}\rangle_{\mathcal{H}}\, db_{n}^{H}(s)\\&=&\sum_{n=1}^{\infty}\sum_{m=1}^{\infty}e_{m}\int_{0}^{t}\langle K_{H}^{*}(f(s)e_{n}), e_{m}\rangle_{\mathcal{H}}\, db_{n}^{H}(s)\\&=&\sum_{n=1}^{\infty}\int_{0}^{t}K_{H}^{*}(f(s)e_{n})\, db_{n}^{H}(s).
        \end{eqnarray*}

\subsection{Smoothing effect}
Consider the linear homogeneous equation associated to \eqref{main}
\begin{equation}\label{HLE}
\partial_tu+(-\Delta)^{\theta/2}u=0,\quad u(0)=u_0.
\end{equation}
It is well known that the operator $(-\Delta)^{\theta/2}$, for $\theta > 0$, generates a semi-group $e^{-t(-\Delta)^{\theta/2}}$ whose kernel $E_{\theta}$ is smooth, radial, and satisfies the scaling property.
\begin{equation}\label{kernel}
E_{\theta}(x,t)=t^{-\frac{N}{\theta}}E_\theta(t^{-\frac{1}{\theta}}x,1):=t^{-\frac{N}{\theta}} K_\theta\left(t^{-\frac{1}{\theta}}x\right),
\end{equation}
where 
\begin{equation}
    \label{K-theta}
    K_\theta(x)=(2\pi)^{-N/2}\int_{\R^N}\,e^{i x\cdot\xi}\,e^{-|\xi|^{\theta}}\,d\xi.
\end{equation}
Hence, the solution of \eqref{HLE} may be formally realized via convolution by 
\begin{equation}
\label{convol}
  u(t,x)=\left(e^{-t(-\Delta)^{\theta/2}}\right)u_0(x)=\bigg(t^{-\frac{N}{\theta}} K_\theta\left(t^{-\frac{1}{\theta}}\cdot\right)\ast u_0\bigg)(x),  
\end{equation}
 whenever this representation makes sense.
 We recall the following point-wise estimate for the kernel $K_\theta$.
 \begin{lemma}{\cite[Lemma 2.1, p. 463]{MYZ2008}}
     \label{K-theta-est}
     Let $\theta>0$ and $N\geq 1$. Then, we have
     \begin{equation}
         \label{kthetaest}
         \left|K_\theta(x)\right|{\leq C}\left(1+|x|\right)^{-N-\theta},\quad x\in\R^N,
     \end{equation}
     {where $C=C(\theta, N)$ is a positive constant.} As a result, we get $K_\theta\in L^r(\R^N)$ for any $1\leq r\leq \infty$.
 \end{lemma}
 \begin{remark}
 {\rm 
    ~\begin{itemize}
         \item[(i)] The inequality \eqref{kthetaest} can be derived from \cite[Theorem 2.1, p. 263]{BG1960} in the case where $\theta < 2$.
    \item[(ii)] For $\theta = 2m$ with $m \geq 1$ being an integer, a more refined estimate than \eqref{kthetaest} is provided in \cite[Proposition 2.1, p. 1325]{GP1}. More precisely, the following estimate holds:
\begin{equation}
\label{kthetaest=2m}
\left|\int_{\R^N}\,e^{i x\cdot\xi}\,e^{-|\xi|^{2m}}\,d\xi\right| \leq C \exp\left(-\kappa_m|x|^{\frac{2m}{2m-1}}\right),
\end{equation}
where 
\begin{equation}
    \label{kappa-m}
    \kappa_m=(2m-1)(2m)^{-\frac{2m}{2m-1}}\;\sin\left(\frac{\pi}{4m-2}\right).
\end{equation}
See also \cite{BD1996, DDY, HWZD}.
\item[(iii)] {The fact that \(K_{\theta}\) belongs to the space $\displaystyle\cap_{1 \leq r \leq \infty}L^{r}(\mathbb{R}^{N})$ implies that the solution to equation \eqref{HLE}, given by \((e^{-t(-\Delta)^{\theta}})u_{0}\), lies in the space \(C([0,\infty), L^{q}(\mathbb{R}^{N})) \cap C((0,\infty), L^{r}(\mathbb{R}^{N}))\) for all \(r \geq q\). This property is crucial, as it allows the use of the space \(C_{T}(L^{r})\) in the proof of Theorem \ref{LWP}.}
    \end{itemize}
    }
 \end{remark}
 From \eqref{convol}, Lemma \ref{K-theta-est} and Young's inequality, we easily derive the following $L^p-L^q$ estimate. See \cite[Proposition 6.1, p. 521]{ADE} for $\theta=4$.
\begin{proposition}\label{prop:smoothing-effect}
Let $\theta>0$ and $N\geq 1$. Then, there exists a positive constant ${\mathcal K}={\mathcal K}(\theta, N)$ such that for all $1\leq p\leq q\leq \infty$, we have
\begin{equation}\label{eq:smoothing-est}
\|e^{-t(-\Delta)^{\theta/2}}\varphi\|_{L^{q}}\leq {\mathcal K}\,t^{-\frac{N}{\theta}(\frac{1}{p}-\frac{1}{q})}\|\varphi\|_{L^{p}},\quad t>0.
\end{equation}
\end{proposition}

 \begin{remark}
     {\rm 
    ~\begin{itemize}
         \item[(i)] By employing the majorizing kernel as derived in \cite[Proposition 2.1, p. 1325]{GP1}, it becomes evident that the constant $\mathcal{H}$ can be chosen independently of both $p$ and $q$ when considering $\theta\in2\N$. This insight has been emphasized and applied in previous studies, including the work presented in \cite{ADE}.
         \item[(ii)] To establish that the constant $\mathcal{H}$ in \eqref{eq:smoothing-est} can be chosen independent with respect to both $p$ and $q$, given any positive $\theta$, we employ an interpolation inequality, allowing us to express it as follows:
         $$
       \|K_\theta\|_{L^r}\leq \|K_\theta\|_{L^\infty}^{1-1/r}\|K_\theta\|_{L^1}^{1/r}\leq \|K_\theta\|_{L^\infty}+\|K_\theta\|_{L^1}:= \mathcal{K}.
         $$
      \item[(iii)] { The semigroup $e^{-t(-\Delta)^{\theta/2}}$ exhibits a contractive property in Lebesgue spaces, which can be understood as follows:
      $$
      \|e^{-t(-\Delta)^{\theta/2}}\|_{L^q-L^q}\leq 1.
      $$}   
     \end{itemize}
      }
 \end{remark}

 For $N\geq 1$, $\theta>0$  and $\gamma\in (0,N)$, define the operator
\begin{equation}
\label{Stheta-alpha}
{\mathbf S}_{\theta,\gamma}(t)=e^{-t(-\Delta)^{\theta/2}}|\cdot|^{-\gamma}.
\end{equation}

\begin{proposition}\label{prop:heat-hardy-kernel-est}
Let $N\geq 1$, $\theta>0$, $\gamma\in (0,N)$ and $1<p,q\leq \infty$  such that
\begin{equation}
\label{eq:relation-p1-p2}
\frac{1}{q}<\frac{\gamma}{N}+\frac{1}{p}<1.
\end{equation}
There exits a positive constant $C:=C(N,p,q,\theta,\gamma)$ such that
\begin{equation}\label{eq:mainest}
\|{\mathbf S}_{\theta,\gamma}(t)\varphi\|_{L^{q}}\leq Ct^{-\frac{N}{\theta}(\frac{1}{p}-\frac{1}{q})-\frac{\gamma}{\theta}}\|\varphi\|_{L^{p}},\quad t>0.
\end{equation}
\end{proposition}
\begin{remark}
     {\rm 
    ~\begin{itemize}
         \item[(i)] The proof of \eqref{eq:mainest} for $\theta=2$ was given in {\cite[Proposition 2.1, p. 121]{BTW}.}
         \item[(ii)] The estimate \eqref{eq:mainest} for $\theta\in 2\N$ was derived in \cite[Proposition 3.8]{G}.
     \end{itemize}
      }
 \end{remark}
\begin{proof}[{Proof of Proposition \ref{prop:heat-hardy-kernel-est}}]
To enhance reader convenience, we provide a comprehensive proof for any $\theta>0$. As one will notice, our proof incorporates arguments borrowed from \cite{BTW, G}. By employing a scaling argument, we can reduce the task of proving \eqref{eq:mainest} to the case of $t=1$. In other words, it is sufficient to establish the inequality
\begin{equation}
\label{eq:mainest=1}
\|{\mathbf S}_{\theta,\gamma}(1)\varphi\|_{L^{q}}\leq C\|\varphi\|_{L^{p}}.
\end{equation}
To achieve this, we employ the following decomposition of the singular weight $|x|^{-\gamma}$:
\begin{equation}
    \label{decomp}
    |x|^{-\gamma}=|x|^{-\gamma} \chi_{\{|x|<1\}}+|x|^{-\gamma} \chi_{\{|x|\geq 1\}}:=f(x)+g(x).
\end{equation}
It is evident that $f\in L^{a}$ holds for all $a<\frac{N}{\gamma}$ and $g\in L^{b}$ holds for all $b>\frac{N}{\gamma}$. Using \eqref{eq:relation-p1-p2}, it becomes apparent that $\frac{N}{\gamma}-\frac{p}{p-1}>0$ and $\frac{\gamma}{N}>\frac{1}{q}-\frac{1}{p}$. Consequently, one can select a sufficiently small  $\varepsilon>0$ such that:
\begin{equation}
    \label{tau-eta}
    \frac{\gamma}{N+\gamma\varepsilon}\geq \frac{1}{q}-\frac{1}{p} \quad\mbox{and}\quad  \frac{\gamma}{N-\gamma\varepsilon}\leq 1-\frac{1}{p}.
\end{equation}
Let us consider $r$ and $s$ with $r, s\geq 1$ satisfying the condition:
\begin{equation}
    \label{r-s}
\frac{1}{r}=\frac{\gamma}{N-\gamma\varepsilon}+\frac{1}{p},\quad
\frac{1}{s}=\frac{\gamma}{N+\gamma\varepsilon}+\frac{1}{p}.
\end{equation}
By employing the smoothing effect \eqref{eq:smoothing-est} and applying H\"older's inequality, we obtain
\begin{equation}
    \label{r-s-est}
    \begin{split}
\|{\mathbf S}_{\theta,\gamma}(1)\varphi\|_{L^q}&\leq \|{\mathbf S}_{\theta, 0}(1)(f\varphi)\|_{L^q}+\|{\mathbf S}_{\theta,0}(1)(g\varphi)\|_{L^q}\\
&{\leq C }\|f\varphi\|_{L^r}+\|g\varphi\|_{L^s}\\
&{\leq C } \left( \|f\|_{L^{\frac{N}{\gamma}-\varepsilon}}+\|g\|_{L^{\frac{N}{\gamma}+\varepsilon}}\right)\, \|\varphi\|_{L^p}.
    \end{split}
\end{equation}
This finishes the proof of Proposition \ref{prop:heat-hardy-kernel-est}.
\end{proof}
\section{Proof of the main result}
\label{S3}
This section is devoted to the proof of Theorem \ref{LWP}.
\begin{lemma}
\label{Parameters}
Let $N, \theta, \gamma, p, q$ be as defined in \eqref{Main-Ass}.  Then, there exists $1 < r < \infty$ such that:
\begin{equation}
\label{r-ass}
\frac{1}{p}\left(\frac{1}{q}-\frac{\gamma}{N}\right)<\frac{1}{r}<\frac{1}{q}.
\end{equation}
Furthermore, the following properties hold 
\begin{eqnarray*}
\label{r-theta}
   && 1-p\sigma>0, \\
    \label{r-theta1}
    &&1-\frac{N}{\theta}\left(\frac{p}{r}-\frac{1}{q}\right)-\frac{\gamma}{\theta}-p\sigma>0,\\
    \label{r-theta2}
    &&1+(1-p)\sigma-\frac{N(p-1)}{r\theta}-\frac{\gamma}{\theta}>0,
\end{eqnarray*}
where
\begin{equation}
\label{sigma}
\sigma = \frac{N}{\theta}\left(\frac{1}{q} - \frac{1}{r}\right).
\end{equation}
\end{lemma}
The proof of Lemma \ref{Parameters} is straightforward. It is worth emphasizing that one can choose
$
    r=\frac{2pq}{p+1}.
$
In addition, with any $r$ satisfying \eqref{r-ass}, we have
\begin{equation}
    \label{r-ineq}
\frac{1}{p}\max\left(\frac{1}{q}-\frac{\gamma}{N},\;\; \frac{p}{q}-\frac{\theta}{N}\right)<\frac{1}{r}<\frac{1}{p}\min\left(1-\frac{\gamma}{N},\;\; \frac{1}{q}-\frac{\gamma}{N}+\frac{\theta}{N},\;\;\frac{p}{q}\right).
\end{equation}
\subsection{Proof of Theorem \ref{LWP}}

First, we consider the linear Cauchy problem
\begin{equation}
\label{Lin}
\left\{
\begin{matrix}
\partial_t{\mathbf Z}(t)+(-\Delta)^{\theta/2} {\mathbf Z}(t)=\mu\,\partial_t B^H(t),\quad t>0,\\
{\mathbf Z}(0)= 0.\\
\end{matrix}
\right.
\end{equation}
The mild solution of \eqref{Lin} is given by
$$
{\mathbf Z}(t)=\mu\int_0^t\,{\rm e}^{-(t-s)(-\Delta)^{\theta/2}}\,d B^H (s).
$$
Since $qH>1,$ $H>\frac{N}{2\theta}$ and $\frac{1}{2}<H<1$, we know from \cite{CMO} that the mild solution ${\mathbf Z}$ belongs to $C([0,T]; L^q)$ for any $T>0$. For $\sigma > 0$ and $1 < r < \infty$, as specified in Lemma \ref{Parameters}, we define
\begin{equation}
\label{K}
{\mathbf K}(T)=\sup_{0\leq t\leq T}\|{\mathbf Z}(t)\|_q+\sup_{0{<} t\leq T}\,\left(t^{\sigma}\,\|{\mathbf Z}(t)\|_r\right).
\end{equation}
Now we are ready to give the detailed proof of Theorem \ref{LWP}. 

We will apply the Banach fixed-point theorem to the integral equation (\ref{integral}) with $\mu\neq 0$. Let us fix $M>2\mathcal{K}\|u_0\|_q$ with $\mathcal{K}$ being given in \eqref{eq:smoothing-est}. For $T>0$, let $\mathbf{X}_T$ be defined as 
\begin{equation}
\mathbf{X}_T=\Bigg\{ u\in C_T(L^q)  \cap C({(0,T]},L^r);\quad  \max\left(\sup_{0\leq t\leq T}\|{u(t)}\|_q , \sup_{0{<}t\leq T}\,\left(t^{\sigma}\,\|{u}(t)\|_r \right)\right)\leq M\Bigg\}.
\end{equation}
We endow $\mathbf{X}_T$ with the metric 
$$
d(u,v)=\sup_{0\leq t\leq T}\|u(t)-v(t)\|_q+\sup_{0{<} t\leq T}\left(t^{\sigma}\|u(t)-v(t)\|_r\right),
$$
where $r, {q}$ {and} $\sigma$ are as in Lemma \ref{Parameters}.
We also define the operator $\Phi:\mathbf{X}_T\to{C_T(L^q)  \cap C({(0,T]},L^r)}$ by
\begin{equation}
\label{Phi}
\Phi(u)(t)=e^{-t(-\Delta)^{\theta/2}}u_0 + \int_0^t {\mathbf S}_{\theta,\gamma} (t-s) (|u(s)|^{p-1}u(s))\, ds+{\mathbf Z}(t),
\end{equation}
where ${\mathbf S}_{\theta,\gamma}$ is given by \eqref{Stheta-alpha}. We aim to demonstrate that $\Phi(\mathbf{X}_T) \subset \mathbf{X}_T$ and that $\Phi$ acts as a contraction under appropriate choice of $T>0$.
First, let us show that $\Phi$ maps $\mathbf{X}_T$ into itself for suitable choice of $T$. For $u\in\mathbf{X}_T$, we have by using \eqref{eq:smoothing-est}, \eqref{eq:mainest} and Lemma \ref{Parameters}
\begin{equation}
    \label{Stab-1}
    \begin{split}
        \|\Phi(u)(t)\|_q&\leq \|e^{-t(-\Delta)^{\theta/2}}u_0\|_q + \|{\mathbf Z}(t)\|_q+\displaystyle\int_0^t \left\|{\mathbf S}_{\theta,\gamma} (t-s) (|u(s)|^{p-1}u(s))\right\|_q\, ds\\
        &\leq \mathcal{K}\|u_0\|_q+\mathbf{K}(T)+C\displaystyle\int_0^t(t-s)^{-\frac{N}{\theta}(\frac{p}{r} - \frac{1}{q})-\frac{\gamma}{\theta}}\,\left\||u(s)|^p\right\|_{\frac{r}{p}}\,ds\\
         &\leq \mathcal{K}\|u_0\|_q+\mathbf{K}(T)+C M^p\displaystyle\int_0^t(t-s)^{-\frac{N}{\theta}(\frac{p}{r} - \frac{1}{q})-\frac{\gamma}{\theta}}\,s^{-p\sigma}\,ds\\
         &\leq \mathcal{K}\|u_0\|_q+\mathbf{K}(T)+C M^p\, T^{a+b-1}\,\mathcal{B}(a,b), 
    \end{split}
\end{equation}
where $\mathcal{B}$ is the standard Beta function and 
\begin{equation}
    \label{a-b}
    a=1-\frac{N}{\theta}\left(\frac{p}{r} - \frac{1}{q}\right)-\frac{\gamma}{\theta},\quad b=1-p\sigma.
\end{equation}
Arguing similarly as above, we obtain
\begin{equation}
    \label{Stab-2}
   t^{\sigma}\|\Phi(u)(t)\|_r \leq \mathcal{K}\|u_0\|_q+\mathbf{K}(T)+C M^p\, T^{a+b-1}\,\mathcal{B}(a-\sigma,b),
\end{equation}
with $a, b$ defined as in \eqref{a-b}. From Lemma \ref{Parameters}, we have $a>0, b>0$, $a+b-1>0$, and $a-\sigma>0$. Therefore, for $T>0$ small enough, {we obtain $\Phi(\mathbf{X}_T) \subset \mathbf{X}_T$ due to the choice of $M > 2\mathcal{K}\|u_0\|_q$.}
Next, we aim to show that $\Phi$ acts as a contraction for $T>0$ small enough. Let $u, v \in \mathbf{X}_T$. Observing that
$$
\left| |u|^{p-1}u-|v|^{p-1}v\right|{\leq p}\,|u-v| \left(|u|^{p-1}+|v|^{p-1}\right),
$$
and using Proposition \ref{prop:heat-hardy-kernel-est} together with H\"older's inequality, we infer
\begin{equation}
    \label{Contra-1}
    \begin{split}
        \|\Phi(u)(t)-\Phi(v)(t)\|_q&\leq\, {C(H)}\displaystyle\int_0^t (t-s)^{a-1-p\sigma}\big\||u(s)|^{p-1}u(s)-|v(s)|^{p-1}v(s)\big\|_{\frac{q}{p}}\,ds\\
        &\leq\, {C(H)}\displaystyle\int_0^t (t-s)^{a-1-p\sigma}\,\|u(s)-v(s)\|_{q}\big(\|u(s)\|_{q}^{p-1}+\|v(s)\|_{q}^{p-1} \big)\,ds\\
        &\leq {C(H)}\ M^{p-1}\,T^{a-p\sigma}\,d(u,v),
    \end{split}  
\end{equation}
where $\sigma$ and $a$  are given by \eqref{sigma} and \eqref{a-b} respectively. Likewise as above, we also obtain
\begin{equation}
    \label{Contra-2}
    t^{\sigma}\|\Phi(u)(t)-\Phi(v)(t)\|_r \leq\,{C(H)}\ M^{p-1}\,T^{a+\sigma}\,d(u,v).
\end{equation}
Plugging estimates \eqref{Contra-1} and \eqref{Contra-2} together, we get
\begin{equation}
    \label{Contraction}
    d(\Phi(u),\Phi(v))\leq\, {C(H)}\ M^{p-1}\left(T^{a-p\sigma}+T^{a+\sigma}\right)\,d(u,v).
\end{equation}
{In the estimates \eqref{Contra-1}, \eqref{Contra-2}, and \eqref{Contraction}, we use the notation $C(H)$ to emphasize the implicit dependence of the constant on the Hurst parameter $H$ associated with the noise $B^H$.}

Given that $a>0$, $\sigma>0$, and $a-p\sigma>0$, we can deduce from \eqref{Contraction} that $\Phi$ exhibits contraction properties for sufficiently small $T>0$. By invoking the Banach fixed point theorem, we successfully conclude the proof of Theorem \ref{LWP}.
\section{{Conclusion and open problems }}

In this paper, we have investigated the local well-posedness of the fractional Hardy-H\'enon equation driven by fractional Brownian motion ({\sf fBm}). By carefully selecting the parameters of the equation and employing the properties of the fractional heat kernel, we established the existence and uniqueness of mild solutions in a suitable Lebesgue spaces. Our results extend previous work on deterministic and stochastic versions of the Hardy-H\'enon equation by incorporating the effects of fractional noise, which introduces additional challenges due to the irregularity of the {\sf fBm}.

The key contributions of this work include the introduction of a mild solution framework for the fractional stochastic Hardy-H\'enon equation, which allows us to handle the singular drift term and the fractional noise simultaneously. We applied the Banach fixed-point theorem to prove local existence and uniqueness under suitable conditions on the parameters of the equation. Additionally, we conducted a careful analysis of the stochastic convolution term, which arises from the {\sf fBm} and plays a crucial role in the well-posedness of the problem.

Despite these advancements, several open problems remain. First, while we have established local well-posedness, the question of global existence and uniqueness of solutions remains open. Extending the local solution to a global one would require additional a priori estimates and a deeper understanding of the long-time behavior of the solution, particularly in the presence of the singular drift term and fractional noise.

Second, the regularity of the solutions, both in space and time, is an important open problem. In particular, it would be interesting to investigate whether the solutions belong to higher-order Sobolev or Besov spaces and how the regularity depends on the Hurst parameter \( H \) of the {\sf fBm}.

Third, for deterministic versions of the Hardy-H\'enon equation, blow-up phenomena have been extensively studied \cite{Ma, Lama, Qi}. It would be valuable to explore whether similar blow-up results hold in the stochastic setting and how the presence of fractional noise influences the blow-up behavior.

Finally, while we have focused on {\sf fBm}, it would be interesting to consider other types of noise, such as L\'evy noise or multiplicative noise. Extending the results to these cases would broaden the applicability of the theory and provide a more comprehensive understanding of stochastic partial differential equations with singular drift terms.

\vspace{0.5cm}

\hrule

\vspace{0.5cm}

\noindent{\bf\large{ Acknowledgements.}} {\em
The authors would like to express their gratitude to the Editor for the meticulous handling of the manuscript and to the Referees for the thorough review and invaluable feedback. The insightful comments and constructive suggestions have greatly improved the quality and clarity of this work.}

\vspace{0.5cm}

\hrule

\vspace{0.5cm}

%%%%%%%%%%%%%%%%%%%%%%%%%%%%%%%%%%%%%%%%%

\end{document}